\input amstex
\documentstyle{amsppt}
\NoRunningHeads
\TagsOnRight

\def\N{\Bbb N}
\def\Z{\Bbb Z}
\def\NC{N_{\Bbb C}}
\def\CH{\chi^\gamma}
\def\i{\bold 1}
\def\Ld{\roman L_2(N',\mu)}
\def\H{\Cal H}
\def\prlim{\mathop{\roman {pr\ lim}}}
\def\indlim{\mathop{\roman {ind\ lim}}}
\def\cg{{\chi,\gamma}}
\def\SQ{$\square$}
\def\l{\langle}
\def\r{\rangle}
\def\L{\langle\!\langle}
\def\R{\rangle\!\rangle_H}
\def\wo{\widehat\otimes}
\def\Ho{\roman{Hol}_0(\NC)}
\def\w{\widetilde}
\topmatter
\title
Dual $n_1$-Appell-like Systems in Infinite-Dimensional Analysis
\endtitle
\author
N. A. Kachanovsky \\
National Technical University of Ukraine "KPI", Kiev, Ukraine
\endauthor
\endtopmatter

\document

The biorthogonal approach to non-Gaussian infinite-dimensional analysis [1]
is generalizing and developing in, at least, three directions. The first
one deals with the biorthogonal decomposition of $\Ld$, where $\mu$ is a
probability measure on co-nuclear space $N'$, by means of dual Appell
systems, which consist of Appell polynomials $\l P_n(x),\varphi^{(n)}\r$,
for which
$$
\int_{N'}\l P_n(x),\varphi^{(n)}\r\mu(\roman dx)=0,\ n\in\N, \tag1
$$
and dual set of distributions (see [2--5]).

The second direction is biorthogonal analysis, connected with hypergroups
([6]) and generalized translation operators ([7]), where the Appell
polynomials are substituted by Appell or Delsarte characters, that,
generally speaking, are not polynomials.

The third direction deals with the generalized Appell ([8,9]) and
Appell-like ([10,11]) systems. These systems do not satisfy
\thetag1 and allow to study pseudodifferential equations on Kondratiev
spaces, see [8--11].

In this paper we will generalize the construction and results of [10,11]
on the wider class of generating functions. The necessity of such
generalization arose because requirements of applications. For example,
this generalization gives the possibility to construct biorthogonal
analysis with respect to measures, which are connected with symmetric
Brownian motions [12].

1. Let $\H$ be a separable real Hilbert space, $N$---separable nuclear
Fr\'echet space, which is embedded in $\H$ topologically (i.e. densely and
continuously). We denote by $N'$ the space, which is dual to $N$ with
respect to $\H$. Further we will use representations
$N=\prlim_{p\in\Bbb N}\H_p$, $N'=\indlim_{p\in\Bbb N}\H_{-p}$ where
${\H}_p$ are Hilbert spaces such that $\forall p$ $\exists p'>p$:
$\H_{p'}\subseteq\H_p$, where embedding is of Hilbert--Schmidt class,
$\H_{-p}$ is dual space of the chain $\H_p\subseteq\H\subseteq\H_{-p}$
(see, for example, [4]). Let us denote by $|\cdot|_p$ the norm in $\H_p$
($p\in\Bbb Z$), $\l x,\theta\r$---dual pairing between $x\in N'$ and
$\theta\in N$, which is given by extension of the inner product in $\H$. We
will preserve these notations for tensor powers of spaces and
complexifications. We will denote by subscript $\Bbb C$ the complexification of
spaces, for example $\NC$ is the complexification of $N$. Let $\Ho$ be the
algebra of germs of functions $\gamma: \NC\to\Bbb C$, holomorphic in $0\in\NC$.
We denote by $\wo$ the symmetric tensor product.

2. Fix $n_1\in\N$. In what follows we denote $\w n:=nn_1$, $n\in\N$. Let
$\chi:\Bbb C\to\Bbb C$ be entire function with $\chi(0)=0$, and its Taylor
decomposition is $\chi(s)=\sum_{n=0}^\infty\frac{\chi_{\w n}}{\w n!}s^{\w
n}$, $s\in\Bbb C,\ \chi_{\w n}\ne0$, $\forall n\in\Z_+$. Let $\gamma\in\Ho$,
$\gamma(0)\ne0$ and in some neighborhood of $0\in\NC$
$\gamma(\theta)=\sum_{n=0}^\infty\frac1{\w n!}\l\gamma_{\w n},
\theta^{\otimes\w n}\r$, $\theta\in\NC$, $\gamma_{\w n}\in\NC^{'\wo\w n}$
(it is not hard to prove, that if $\gamma^{(1)}$ and $\gamma^{(2)}$ satisfy
these conditions, then $\gamma^{(1)}(\theta)/\gamma^{(2)}(\theta)$ also
satisfy them). We put $\CH(\theta;z):=\gamma(\theta)\chi(\l z,\theta\r)$,
$z\in\NC'$, $\theta\in\NC$. One can decompose $\CH(\cdot;z)$ in Taylor
series and obtain by kernel theorem the representation
$$
\CH(\theta;z)=\sum_{n=0}^\infty\frac1{\w n!}
\l P^\cg_{\w n}(z),\theta^{\otimes\w n}\r,
\ P^\cg_{\w n}(z)\in\NC^{'\wo\w n}.
$$
\definition{Definition 1}
We will call $\bold P^\cg:=\big\{\l P^\cg_{\w n}(\cdot),\varphi^{(\w n)}\r:$
$\varphi^{(\w n)}\in\NC^{\wo\w n}$, $n\in\Z_+\big\}$ the (generalized)
$n_1$-Appell-like polynomial system.
\enddefinition

Particularly, $1$-Appell-like polynomials are Appell-like polynomials, see
[10,11] and references therein.

The properties of $n_1$-Appell-like polynomials, as it easy can be proved,
are analogous to the same of Appell-like polynomials. So, for
example,
$$
P^{\chi,\gamma^{(1)}}_{\w n}(z)=\sum_{m=0}^nC_{\w n}^{\w m}
P^{\chi,\gamma^{(2)}}_{\w m}(z)\wo\widehat\gamma_{\w{n-m}}, \tag2
$$
where $\gamma^{(1)}$, $\gamma^{(2)}$ satisfy the conditions above to
$\gamma$, $\widehat\gamma_{\w{n-m}}\in\NC^{'\wo (\w{n-m})}$ from
decomposition
$\gamma^{(1)}(\theta)/\gamma^{(2)}(\theta)=\sum_{n=0}^\infty\frac1{\w n!}
\l\widehat\gamma_{\w n},\theta^{\otimes\w n}\r$.

Let $\Cal P_\cg:=\big\{\varphi(x)=\sum_{n=0}^{M_\varphi}
\l P^\cg_{\w n}(x),\varphi^{(\w n)}\r$: $x\in N'$,
$\varphi^{(\w n)}\in\NC^{\wo\w n}\big\}$. It follows from \thetag2 that
$\Cal P_\cg$ does not depend on $\gamma$. The dependence on $\chi$
connected only with $n_1$, therefore for $\Cal P_\cg$ we can accept the
notation $\Cal P_{n_1}$. We choose the topology in $\Cal P_{n_1}$ such that
it be topologically isomorphic to topological direct sum of tensor powers
of $\NC$ $\underset{n=0}\to{\overset\infty\to\oplus}\NC^{\wo\w n}$ via the
coefficients of $n_1$-Appell-like polynomials. Note, that if $n_1\ne1$, then
$\Cal P_{n_1}$ does not coincide with the set $\Cal P=\Cal P_1$
of all continuous polynomials on $N'$, because $\Cal P_{n_1}$ contains only
polynomials with $n_1$-divisible powers.

Let $(\H_p)^1_{q,\cg}$ be the completion of $\Cal P_{n_1}$ with respect to
the norm $\|\cdot\|_{p,q,\cg}$, which is defined for
$\varphi(x)=\sum_{n=0}^\infty\l P^\cg_{\w n}(x),\varphi^{(\w n)}\r$, $x\in
N'$, $\varphi^{(\w n)}\in\NC^{\wo\w n}$ by equality $\|\varphi\|_{p,q,\cg}^2
:=\sum_{n=0}^\infty(\w n!)^22^{q\w n}|\varphi^{(\w n)}|_p^2$;
$(N)^1_\cg:=\prlim_{p,q\in\N}(\H_p)^1_{q,\cg}$.

\proclaim{Theorem 1} For $\cg$ under conditions above the
spaces $(N)^1_\cg$ are invariant with respect to $\gamma$,
i.e. $(N)^1_\cg=(N)^1_\chi$. If there exists positive nondecreasing
function $\Psi$ of real argument such that
$|\chi(z)|\leqslant\Psi(|z|)$, $\forall z\in\Bbb C$, then each element
from $(N)^1_\chi$ is restriction on $N'$ of some entire on $\NC'$ function.
Moreover, if there exist $C_0>0$ and $C_1>0$ such that $\Psi(r)\leqslant
C_0\exp\{C_1r\}$, then $(N)^1_\chi\subseteq(N)^1:=(N)^1_{\exp}$.
\endproclaim
\noindent
The proof of this theorem is based on formula \thetag2, for details see
[11]. \SQ

\example{Example}
Let us consider the function $\Lambda_s(u)
:=\sum_{m=0}^\infty\frac{u^{2m}}{(2m)!}\cdot\big(-\frac14\big)^m
\frac{(2m)!s!}{m!(s+m)!}$, $u\in\Bbb C$, $s\in[-\frac12,\infty)$, where
$s!:=\Gamma(s+1)$ ($\Lambda_s$ is related to random walks with spherical
symmetry, see [12]). It is easy to prove, that $\Lambda_s(u)$ satisfies all
conditions of Theorem 1: $n_1=2$, $\Psi(r)=e^r$.
\endexample

We define on $\Cal P_{n_1}$ pseudodifferential operator
$\l\varPhi^{(\w n)},D_\chi^{\otimes\w n}\r$, $\varPhi^{(\w n)}\in
\NC^{'\wo\w n}$, setting on monomials
$$
\multline
\l\varPhi^{(\w n)},D_\chi^{\otimes\w n}\r\l x^{\otimes\w m},
\varphi^{(\w m)}\r\\
:=1_{\{m\geqslant n\}}\frac{\w m!\chi_{\w{m-n}}}{(\w{m-n})!\chi_{\w m}}
\l x^{\otimes (\w{m-n})}\wo\varPhi^{(\w n)},\varphi^{(\w m)}\r,
\ \varphi^{(\w m)}\in\NC^{\wo\w m},\ x\in N',
\endmultline
$$
and extending by linearity (here $1_{\{m\geqslant n\}}$ is the indicator of
$\{m\geqslant n\}$).

\proclaim{Theorem 2}
The $n_1$-Appell-like polynomials are generalized powers with respect to
$\l\varPhi^{(\w n)},D_\chi^{\otimes\w n}\r$, i.e.
$$
\l\varPhi^{(\w n)},D_\chi^{\otimes\w n}\r
\l P^\cg_{\w m}(x),\varphi^{(\w m)}\r
=1_{\{m\geqslant n\}}\frac{\w m!}{(\w{m-n})!}\l
P^\cg_{\w{m-n}}(x)\wo\varPhi^{(\w n)},\varphi^{(\w m)}\r.
$$
\endproclaim
\noindent
The proof uses \thetag2 and definition of
$\l\varPhi^{(\w n)},D_\chi^{\otimes\w n}\r$. \SQ

\proclaim{Corollary}
For all $\varPhi^{(\w n)}\in\NC^{'\wo\w n}$
$\l\varPhi^{(\w n)},D_\chi^{\otimes\w n}\r\in\Cal L(\Cal P_{n_1},
\Cal P_{n_1})$.
\endproclaim

3. Let $H$ be a separable real Hilbert space such that $\Cal P_{n_1}$
is embedded in $H$ and for $\varphi\in\Cal P_{n_1}$
$\|\varphi\|_H=0\Rightarrow\varphi(x)=0\ \forall x\in N'$. Suppose also
that there exist $C>0$, $K>0$, $p\in\N$ such that
$\||P^{\chi,\i}_{\w n}(\cdot)|_{-p}\|_H\leqslant\w n!C^{\w n}K$, and embedding
$\Cal P_{n_1}\hookrightarrow H$ is topological.

\remark{Remark}
Note, that in the case $n_1=1$ the main example of $H$ is $\Ld$, where
$\mu$ is analytic nondegenerate probability measure (see [4,5,8--11]). But
in the case $n_1\ne1$ such spaces can turn out too wide for $H$.
For example, if $\mu$---smooth measure with second analyticity
condition (see [2]), then $\Cal P_{n_1}$ is not dense in $\Ld$, if
$n_1\ne1$. But in the similar cases one can use the completion of
$\Cal P_{n_1}$ with respect to the norm of $\Ld$ as $H$.
\endremark

Let $\Cal P_{n_1,H}'$ be the space, which is dual to $\Cal P_{n_1}$ with
respect to $H$. Let $\l\varPhi^{(\w n)},D_\chi^{\otimes\w n}\r^*_H$
$\in\Cal L(\Cal P_{n_1,H}',\Cal P_{n_1,H}')$ be the operator, which is dual
to $\l\varPhi^{(\w n)},D_\chi^{\otimes\w n}\r$, i.e. for all $\Phi\in
\Cal P_{n_1,H}'$, $\varphi\in\Cal P_{n_1}$
$\L\l\varPhi^{(\w n)},D_\chi^{\otimes\w n}\r^*_H\Phi,\varphi\R
=\L\Phi,\l\varPhi^{(\w n)},D_\chi^{\otimes\w n}\r\varphi\R$, where
$\L\cdot,\cdot\R$ is the dual pairing between $\Cal P_{n_1,H}'$ and
$\Cal P_{n_1}$, which is given by extension of the inner product in $H$.
One can prove that there exist $p'=p'(\cg)$, $q'=q'(\cg)$ such that for all
$p\geqslant p'$, $q\geqslant q'$ $(\H_p)^1_{q,\cg}\hookrightarrow H$ 
topologically.
Therefore, $(N)^1_\chi\hookrightarrow H$ topologically. Let
$(N)^{-1}_{\chi,H}$ be the space, which is dual to $(N)^1_\chi$ with
respect to $H$. We define on $(N)^1_\chi$ the linear functional $\delta_z$,
$z\in\NC'$, by formula $\L\delta_z,\varphi\R:=\varphi(z)$,
$\varphi\in(N)^1_\chi$. It is not hard to prove that $\delta_z\in
(N)^{-1}_{\chi,H}$. Let $\gamma$ be such as above. We put
$$
Q^\cg_{H,\w m}(\varPhi^{(\w m)};\cdot):=\sum_{k=0}^\infty\frac1{\w k!}
(\l\varPhi^{(\w m)}\wo\w\gamma_{\w k},D_\chi^{\otimes (\w{m+k})}\r^*_H
\delta_0)(\cdot)\in\Cal P_{n_1,H}', \tag3
$$
where $\varPhi^{(\w m)}\in\NC^{'\wo\w m}$, $\w\gamma_{\w k}\in\NC^{'\wo\w
k}$ from decomposition $1/\gamma(\theta)=\sum_{k=0}^\infty\frac1{\w k!}
\l\w\gamma_{\w k},\theta^{\otimes\w k}\r$, $\theta\in\NC$.

\proclaim{Theorem 3 {\it (biorthogonality w.r.t. $H$)}}
The generalized functions $Q^\cg_{H,\w m}(\varPhi^{(\w m)};\cdot)$ are
biorthogonal to $\bold P^\cg$-system, i.e.
$$
\L Q^\cg_{H,\w m}(\varPhi^{(\w m)};\cdot),
\l P^\cg_{\w n}(\cdot),\varphi^{(\w n)}\r\R
=\delta_{mn}\w n!\l\varPhi^{(\w n)},\varphi^{(\w n)}\r,
\ \varPhi^{(\w m)}\in\NC^{'\wo\w m},\ \varphi^{(\w n)}\in\NC^{\wo\w n}.
$$
\endproclaim
\noindent
The proof is based on \thetag2, \thetag3 and definition of $\delta_0$. \SQ

\proclaim{Theorem 4}
Each generalized function $\Phi\in\Cal P_{n_1,H}'$ can be represented in
the form of
$$
\Phi(\cdot)=\sum_{m=0}^\infty Q^\cg_{H,\w m}(\varPhi^{(\w m)};\cdot),
\ \varPhi^{(\w m)}\in\NC^{'\wo\w m}, \tag4
$$
where the sequence of kernels $\big\{\varPhi^{(\w m)}\big\}_{m=0}^\infty$
is uniquely determined by $\Phi$. Conversely, each sequence
$\big\{\varPhi^{(\w m)}\big\}_{m=0}^\infty$,
$\varPhi^{(\w m)}\in\NC^{'\wo\w m}$ defines the generalized function
$\Phi\in\Cal P_{n_1,H}'$ by formula \thetag4.
\endproclaim
\noindent
The proof is the same as the proof of corresponding statement in [11]. \SQ

\definition{Definition 2}
The system of $n_1$-Appell-like polynomials $\bold P^\cg$ and the family of
generalized functions $Q^\cg_{H,\w m}(\varPhi^{(\w m)};\cdot)$,
$\varPhi^{(\w m)}\in\NC^{'\wo\w m}$, form dual $n_1$-Appell-like system.
\enddefinition

Let us introduce an analog of $S$-transform (see [2--11]). For
$\Phi\in(N)^{-1}_{\chi,H}$ we put
$$
(S_{\cg,H}\Phi)(\theta):=\L\Phi(\cdot),\CH(\theta;\cdot)\R,
\ \theta\in\NC.
$$
One can prove, that this definition is correct.

\proclaim{Theorem 5 {\it (characterization theorem)}}
$S_{\cg,H}$-transform is topological isomorphism from $(N)^{-1}_{\chi,H}$
to subspace of $\Ho$, which consist of functions $F$, which can be written
in the form of $F(\theta)=\sum_{n=0}^\infty\l\varPhi^{(\w n)},
\theta^{\otimes\w n}\r$, $\varPhi^{(\w n)}\in\NC^{'\wo\w n}$.
\endproclaim
\noindent
The proof is analogous to the proof of characterization theorem
in [11].  \SQ

\bigskip
{\smc N. A. Kachanovsky, postgraduate student, National Technical
University of Uk\-ra\-i\-ne "Kiev Polytechnic Institute", Chair of
Mathematical Methods of System Analysis, pr. Pobedy, 37, Kiev, Ukraine.
}

\enddocument